\documentclass[12pt]{article}
\usepackage{cmap}
\usepackage[T2A]{fontenc}
\usepackage[cp866]{inputenc}
\usepackage[english]{babel}
\usepackage[tbtags]{amsmath}
\usepackage{amsfonts,amssymb}

\newcommand\np{\noindent}

\newcommand\CC{\mathbb C}
\newcommand\DD{\mathbb D}

\newcommand\NN{\mathbb N}

\newcommand\ZZ{\mathbb Z}
\newcommand\RR{\mathbb R}
\newcommand\SA{\mathbb S}
\newcommand\GA{\Gamma}
\newcommand\ga{\gamma}
\newcommand\sub{\subset}
\newcommand\sm{\setminus}
\newcommand\pa{\partial}

\newcommand\la{\lambda}
\newcommand\La{\Lambda}
\newcommand\vp{\varphi}
\newcommand\eps{\varepsilon}

\newcommand\HC{\mathcal H}

\newcommand\SI{\Sigma}

\newcommand\trl{\triangleleft}
\newcommand\trr{\hfill\square}

\renewcommand{\abstract}{\small}

\begin{document}
\title{On the removable singularities\\
of complex analytic sets}
\author{E.M.Chirka
\thanks{The work is supported by the Russian Science Foundation
under grant 14-05-00005.}\\
Steklov Math.Institute RAS\\
chirka@mi.ras.ru}
\date{}
\maketitle

\begin{abstract}
There is proved the sufficiency of several conditions
for the removability
of singularities of complex-analytic sets in domains of $\,\CC^n.$
\end{abstract}
\vskip4mm

{\bf 1. Introduction.} A closed subset $\,\SI\,$ of a complex
manifold $\,M\,$ is called below {\sf $p$-removable\/} if for every
respectively closed purely $p$-dimensional complex analytic subset
$\,A\sub M\sm\SI\,$ its closure in $\,M\,$ is an analytic set.
We omit "complex"$ $ in what follows.
Thus "analytic set"$ $ below
means a set $\,A\sub M\,$ such that for every point $\,a\in A\,$
there exists a neighborhood $\,U\sub M\,$ of $\,a\,$ such that
$\,A\cap U\,$ is the set of common zeros of a family of holomorphic
functions in $\,U$.
Such a set $\,A\,$ has in general singular
points but the set $\,sng\,A$ of them is removable in sense of
given definition and in essence we study below the boundary sets
$\,\Sigma\,$ the adding of which to $\,A$ does not spoil the
analyticity.

By $\,\HC^m$ we denote the Hausdorff measure of dimension $\,m\ge
0$ (see e.g. \cite{G} Ch.III).
There is well-known sufficient
metrical condition (Shiffman theorem): {\it $\,\SI\,$ is
$p$-removable if $\,\HC^{2p-1}(\SI)=0$}; see \cite{Sf} or \cite{C}
\S$\,$4.4.
We assume some smooth metric on $\,M\,$ being fixed and
Hausdorff measures are taken with respect to this metric.
The
vanishing of $\,\HC^m$-measure of a set in $\,M\,$ does not depend
on the choice of smooth metric.

The condition
of Shiffman theorem is non-improvable in the scale of Hausdorff
measures and for its weakening one needs in additional assumptions.
For $\,p=1\,$ the following remarkable theorem is obtained by Lee
Stout (\cite{S} Theorem 3.8.18): \vskip2mm

\np{\it Let $\,D\,$ be an arbitrary domain in $\,\CC^n$ and $\,A\,$
is relatively closed purely one-dimensional analytic subset of
$\,D\sm E\,$ where $\,E\,$ is a compact in $\,\bar D\,$ such that
$\,\HC^2(E)=0,\ H^1(E,\ZZ)=0\,$ and the set $\,E\cap\pa D\,$ is
either empty or a point.
Then $\,\bar A\cap D\,$ is one-dimensional
analytic set.} \vskip1mm

This is the first statement on removable singularities which I know
with conditions on the boundaries of tested sets.

Here $\,H^1$ denotes $\check{\rm C}$ech cohomology and the
condition $\,H^1(E,\ZZ)=0\,$ is purely topological.
By Bruschlinsky
theorem it is equivalent to the condition that any continuous
function without zeros on $\,E\,$ has continuous logarithm (see
\cite{B} or \cite{S} p.19).
Such sets are called also {\sf simply co-connected}.
There are for instance all totally disconnected
compact sets (the connected components are points), simple Jordan
arcs, plane compact sets with connected complements e.t.c.

Some conditions on $\,\bar\SI\cap\pa D\,$ are necessary in general
for the removability of $\,\SI\sub D$.
If for example $\,D\,$ is
the unit ball in $\,\CC^n$ and $\,\SI\,$ is its diameter on
the axis $\,x_1=Re\,z_1\,$ then $\,\bar\SI\cap\pa D\,$ consists of two
points only but for the semi-disk $\,A=D\cap\{Im\,z_1>0,
z_2=...=z_n=0\}$ the set $\,\bar A\cap D\,$ is not analytic.
Nevertheless, the condition $\,\#\,(E\cap\pa D)\le 1\,$ in Stout theorem can be
weakened in the following way.
\vskip2mm

{\bf Theorem 1.} {\it Let $\,D\,$ be an arbitrary domain in
$\,\CC^n$ and $\,\SI\,$ be its bounded relatively closed subset.
Assume that $\HC^2(\bar\SI)=0\,$ and the one-point
compactifi\-cation $\,\SI\sqcup\circ\,$ is simply co-connected.
Then $\,\SI\,$ is $1$-removable.} \vskip2mm

One-point compactification of a topological space $\,X\,$
is the topological space $\,X\sqcup\,\circ\,$ with the same topology
on $\,X\,$ and additional point
$\,\circ\,$ which punctured neighbor\-hoods are the complements to
compact subsets of the space $\,X$.
(Here and everywhere below
$\,\sqcup{}\,$ means the union of disjoint sets.)
If $\,\bar\SI\sm\SI\,$ is a point (as in Stout theorem) then
$\,\SI\sqcup\circ=\bar\SI$ as topological spaces, and in the
example with diameter of a ball the set $\,\SI\sqcup\circ\,$
is homeomorphic to a circle and thus is not simply co-connected.
The condition $\,H^1(X\sqcup\circ\,,\ZZ)=0\,$ for relatively closed
$\,X\sub D\Subset\RR^N$ is evidently equivalent to the condition
that every continuous function without zeros on $\,\bar X\,$ and
constant on $\,\bar X\sm X\,$ has on $\,\bar X\,$ a continuous
logarithm which is also constant on $\,\bar X\sm X$. \vskip3mm

The problems with singularities for analytic sets of dimension
$\,p>1\,$ are in some sense simpler due to Harvey -- Lawson theorem
(see for instance theorem 20.2 on compact singularities in \cite{C}),
but the proofs need in additional pseudo\-convexity type assumptions.
The following theorem is valid for arbitrary $\,p\ge 1$.\vskip2mm

{\bf Theorem 2.} {\it Let $\,D\,$ be a domain in $\,\CC^{n},\
\SI\,$ be its relatively closed bounded subset,
$\SI_{b}:=\bar\SI\cap\pa D$ and $\,p\in\NN,\, p\le n$.
Assume that $\,\HC^{2p}(\SI)=0,\
H^{2p-1}(\SI\sqcup\circ,\ZZ)=0\,$ and
$\,\widehat{\SI_b}\cap\SI=\varnothing\,$. Then $\,\SI\,$ is
$p$-removable.}\vskip2mm

Here $\,\hat X\,$ means the polynomially convex hull of a compact set
$\,X\sub\CC^n$ that is the set
$\{z:|P(z)|\le max_X\,|P|\ \text{for any polynomial}\ P\,\}$.
\vskip2mm

Remark that the Theorems, Propositions and Corollaries below are
extending in obvious way onto domains in Stein manifolds (instead
of $\,D\sub\CC^n$ as in the text) in view of the proper
imbeddebility of such manifolds into suitable $\,\CC^{\,N}$.
One should only to substitute the polynomials and polynomially (or
rationally) convex hulls by global holomorphic (meromorphic)
functions and corresponding hulls.
If $\,M\,$ is a complex
submanifold in $\,\CC^{\,N}$, $G$ is a domain on $\,M\,$ and
$\,\SI\,,\,A\sub G$ then $\,D=\CC^{\,N}\sm(M\sm G)$ is a domain in
$\,\CC^{\,N}$, the sets $\SI,\,A\,$ are contained in $\,D\,$ and
one can apply the results in $\,\CC^{N}$ restricting to $\,M$.

The paper is organized as follows.
In sect.2 we prove mainly topological preliminaries in spirit of
argument principle and degrees of continuous mappings.
Sect.3 contains the proof of Theorem 1 in more general situation
of Proposition 1.
Examples in sect.4 stress that simple sufficient conditions of
Theorem 1 are not necessary at all.
The proof of Theorem 2 is placed in sect.5 and the corollaries of
both theorems are collected in sect.6.
At the end we discuss natural relations with removable singularities
of holomorphic (and meromor\-phic) functions and stste some open
questions.
\vskip1mm

There are many references in the text on the book of E.L.Stout
\cite{S} but the given  proofs are complete; the references
indicate simply the corresponding statements and arguments in
\cite{S}. \vskip3mm

{\bf 2. Preliminaries.} First of all we endow
$\,\Sigma\sqcup\circ\,$ by the structure of metrical space inducing
the same topology as described above.
Let $\,\Sigma\,$ be relatively closed and bounded subset of a domain
$\,D\sub\CC^n$.
The point $\,\circ\,$ is connected and we have to
transform $\,\Sigma_b$ into a connected set.
Define $\,\Sigma_0:=\{(tz,t):z\in\Sigma_b\,,\,0\le t\le 1\}$.
Let $\,\vp\,$ be a continuous nonnegative function in
$\,\CC^n\times\RR_{\,t}$ with zero set $\,\Sigma_0$ and
$\,\tilde\rho(x,x'):=\inf\int_\ga\vp\,ds\,$ where $\,ds\,$ is the
euclidean metric in $\,\CC^n\times\RR\,$ and the infimum is taken
by all smooth curves $\,\ga\,$ containing the points $\,x,x'$.
Then $\,\tilde\rho\,$ is symmetric and satisfies the triangle
inequality.
It is degenerated on $\,\Sigma_0\times\Sigma_0$ but the
corresponding distance function
$\,\rho(z,z'):=\tilde\rho((z,1),(z',1))$ for $\,z,z'\in\Sigma,\
\rho(z,\circ):=\tilde\rho((z,1),\Sigma_0)$ for $\,z\in\Sigma\,$ and
$\,\tilde\rho(\circ,\circ):=0\,$ defines a metric on
$\,\Sigma\sqcup\circ\,$ as we need. \vskip3mm

The following reduction is used in both proofs. \vskip2mm

{\bf Lemma 1.} {\it Let $\,\Sigma\,$ be relatively
closed subset of a domain
$\,D\sub\CC^n$ such that $\HC^{2p}(\Sigma)=0\,$ and
$\,A\sub D\sm\Sigma\,$ is relatively
closed purely $p$-dimensional analytic
subset such that $\,\overline{A_0}\cap D\,$ is analytic for any
irreducible component $\,A_0$ of $\,A$.
Then $\,\bar A\cap D\,$ is analytic.}
\vskip2mm

$\trl$ Let $\,A=\cup A_j$ be the decomposition onto
irreducible components and $\,a\in\bar A\cap\Sigma$.
As
$\,\HC^{2p}(\Sigma)=0\,$ there is a complex plane $\,L\ni a\,$ of
complex dimension $\,n-p\,$ such that the set $\,\Sigma\cap L\,$ is
locally finite.
Without loss of generality we can assume that
$\,a=0\,$ and $\,L\,$ is the coordinate plane
$\,z':=(z_1,...,z_p)=0$. Then there is $\,r>0\,$ such that
$\,\Sigma\cap L\cap\{|z|\le r\}=\{0\}$ and $\,L\cap\{|z|\le r\}\sub
D$.
Let us show that there exists a neighborhood $\,U\ni a\,$
intersecting only finite number of $\,A_j$.

Assume not.
Then there is a sequence of points $\,a_k\to a\,$ such
that $\,a_k\in A_k$ and $\,A_k\not=A_l$ for $k\not=l$.
As the
decomposition $\,A=\cup A_j$ is locally finite in $\,D\sm\Sigma\,$
we can assume (passing to a subsequence) that there is $\,r'>0\,$
such that $\,A_k\cap\{|z'|\le r',\,|z|=r\}=\varnothing$ for all
$\,A_k\ni a_k$.
Then the restrictions of the projection $\,z\mapsto
z'$ onto $\,A_k\cap\{|z'|<r',|z|<r\}$ are proper, in particular,
their images contain the ball $\{|z'|<r'\}\sub\CC^{\,p}$.
The projection of
$\,\Sigma\,$ into $\,\CC^{\,p}_{z'}$ has zero volume and thus there
is a point $\,b'\not\in z'(\Sigma)$ with $|b'|<r'$.
By the
construction there are points $\,b_k\in A_k$ such that
$\,z'(b_k)=b'$.
Passing to a subsequence we can assume that there
is $\,b\in D\,$ such that $\,b_k\to b\,$ as $\,k\to\infty$.
But
then $\,b\in D\sm\Sigma\,$ and we obtain the contradiction with
local finiteness of the decomposition into irreducible components
(in a neighborhood of $\,b$).

Thus there is a neighborhood $\,U\ni a\,$ in $\,D\,$ and a finite
number of indexes $\,j_1,...,j_N$ such that $\,U\cap
A_j=\varnothing\,$ if $\,j\not\in\{j_1,...,j_N\}$.
By the condition
the sets $\,\overline{A_{j_\nu}}\cap D\,$ are analytic and thus
$\,\bar A\cap U=\cup_{\nu=1}^N(\overline{A_{j_\nu}}\cap U)$ is also
analytic.
As $\,a\in\bar A\cap\Sigma\,$ is arbitrary the set
$\,\bar A\cap D\,$ is analytic. $\trr$ \vskip3mm

For the proofs of main results we need in the
following lemmas in a spirit of argument principle. \vskip2mm

{\bf Lemma 2.} {\it Let $\,E\,$ be a compact subset of
zero $\,\HC^{m}$-measure in the closure of a
domain $\,D\sub\RR^{N},\,N>m\in\NN,$ and $\,K\sub E\,$ is compact.
Then every continuous map $\,f:K\to\RR^m\sm 0\,$
which is constant on $\,K\cap\pa D\,$
extends to a continuous map of $\,E\to\RR^m\sm 0\,$ which is
constant on $\,E\cap\pa D$.} \vskip2mm

$\trl$ (see \cite{S} Lemma 3.8.16). As $\,0\not\in f(K)$ there is
constant $\,\delta>0,$ such that $|f|>\delta\,$ on $\,K$.
Then
there exists a continuous map $\,\tilde f:\RR^N\to\RR^m\,$
which is equal to $\,f\,$ on $\,K,$ to a constant $\,a\,$ on
$\,E\cap\pa D,\,|a|>\delta,$ and smooth on the set $\{|\tilde
f|<\delta/2\}$.

The image of a set with zero $\,\HC^m$-measure by a smooth
map to $\,\RR^m\,$ also
has zero $\,\HC^m$-measure, hence there is $\,c\in\RR^m,\ |c|<\delta/2\,$, such
that the map $\,\tilde f-c\,$ does not have zeros on $\,E$.
Let $\,\vp:\RR^m\to\RR^m$ be a diffeomorphism
such that $\,\vp(x)=x\,$ if $|x|>\delta\,$ and $\,\vp(c)=0$.
Then $\,\vp\circ\tilde f\,$ is a desired extension.
$\trr$ \vskip3mm

{\bf Corollary 1.} {\it If $\,H^{m-1}((E\cap D)\sqcup\circ,\ZZ)=0\,$ then
$\,H^{m-1}(X\sqcup\circ,\ZZ)=0\,$ for any relatively closed
subset $\,X\sub E\cap D$.}
\vskip2mm

$\trl$ The notion of Hausdorff measures is well defined for an
arbitrary metric space.
With the metric on $\,(E\cap
D)\sqcup\,\circ\,$ defined above we have evidently $\,\HC^m((E\cap
D)\sqcup\,\circ)=\HC^m(E)=0$, hence $\,\HC^m(X\sqcup\,\circ)=0$.
With this property, the condition
$\,H^{m-1}(X\sqcup\,\circ,\ZZ)=0\,$ is equivalent to that every
continuous map $\,f\,$ of $\,X\sqcup\,\circ\,$ into the unit sphere
$\,\SA^{m-1}$ in $\,\RR^m$ is homotopic to a constant one in the
class of continuous mappings into $\,\SA^{m-1}$ (see Theorems VII.3
and VIII.2 in \cite{HW}).
The "projection"$ $ $\,\pi:\bar X\to
X\sqcup\circ\,$  such that $\,\pi(x)=x\,$ for $\,x\in X\,$ and
$\,\pi(x)=\circ\,$ for $\,x\in\bar X\sm X\,$ is continuous.
The map
$\,f\circ\pi:\bar X\to\SA^{m-1}$ is constant on $\,\bar X\sm X$.
By
Lemma 2 there is a continuous map $\,\tilde f: E\to\RR^m\sm 0\,$
equal to $\,f\circ\pi\,$ on $\,\bar X\,$ and to a constant
(${}\equiv\,a$) on $\,E\cap\pa D$.
Set $\,F:=\tilde f/|\tilde f|$ on $\,E\cap D\,$
and $\,F(\circ):=a/|a|$.
As $\,H^1((E\cap D)\sqcup\circ,\ZZ)=0\,$
the map $\,F$ is homotopic to a constant one in the class of
continuous mappings into $\,\SA^{m-1}$.
The same is true for
$\,F|_{X\sqcup\,\circ}$ and thus $\,H^{m-1}(X\sqcup\,\circ,\ZZ)=0$.
$\trr$ \vskip3mm

{\bf Lemma 3.} {\it Let $\,E'\sub E\,$ be compact sets in $\,\CC^n$ and
$\,A\sub\CC^n\sm E\,$
is a bounded purely $\,p$-dimensional analytic set, $1\le p\le n$, with
$\,\pa A\sub E$.
Assume that

$1)\ \HC^{2p}(E\sm E')=0\,$ and

$2)\ \HC^{2p-1}(E')=0\,$ or $\,A\not\sub\widehat{E'}$. \vskip1mm

\np Then $\,H^{2p-1}((E\sm E')\sqcup\circ,\ZZ)\not=0$.} \vskip2mm

$\trl$ (See Theorem 3.8.15 in \cite{S}.) Consider first the case
when $\,\HC^{2p-1}(E')=0$.
Then $\,\HC^{2p}(E)=0\,$ and there is an
affine map $\,f:\CC^{n}_{z}\to\CC^{p}_{w}$ such that $\,0\in f(A)$
but $|f|>r>0\,$ on $\,f(E)$.
As $\,\HC^{2p-1}(E')=0\,$ then one can
assume also that the ray $\,Im\,w_1=0,\,Re\,w_1<0,\,w_j=0,\,1<j\le p,$
does not intersect $\,f(E')$. \vskip1mm

\np Then there exists a homotopy $\,\vp_t:\CC^p\to\CC^p,\,0\le t\le
1,$ such that $\,\vp_0(w)\equiv w,\,\vp_t(w)=w\,$ if $|w|<r,\
|\,\vp_t(w)|\ge r\,$ if $|w|\ge r\,$ and $\,\vp_1\equiv 1\,$ on
$\,f(E')$.
Set $\,F:=\vp_1\circ f$.
Analytic set $\,A\cap\{F=0\}$
is compact and so is finite.
Hence there exists $\,\eps\in(0,r)\,$
and $\,k\in\NN\,$ such that $\,F: A\cap\{|F|<\eps\}\to\{|w|<\eps\}$
is $\,k$-sheeted analytic covering (see \cite{C}). \vskip1mm

Let $\,\rho\,$ be a smooth real function in $\,\CC^{n}$ equal to
$|F|$ when $|F|<\eps,\ =1\,$ on $\,E\,$ and $\,0\le\rho <1\,$ on
$\,A$.
By Sard theorem for analytic sets (Proposition 14.3.1 in
\cite{C}) $A_{t}:=A\cap\{\rho<t\}$ is analytic set with a border
for almost all $\,t\in(0,1)$.
Let $\,t_{j}, j=0,1,...,$ be increasing
sequence of such values with $\,t_j\to 1\,$ as $\,j\to\infty,
\ t_{0}<\eps$ and $\,\GA_{j}:=A\cap\{\rho=t_{j}\}$.
\vskip1mm

Let $\,\theta:=d^{c}\log|F|^{2}\wedge(dd^{c}\log|F|^{2})^{p-1}$
where $\,d^{c}:=i(\bar\pa-\pa)$.
As $\,(dd^{c}\log|w|^{2})^{p}=0\,$
in $\,\CC^{p}_w\sm 0\,$ we have $\,d\theta=0\,$ on $\,A\sm\{F=0\}$.
By Stokes theorem for analytic sets (Theorem 14.3 in \cite{C})
$$
\int_{\GA_{j}}\theta=\int_{\GA_{0}}\theta=
t_{0}^{-2p}\int_{\GA_{0}} d^{c}|F|^{2}\wedge(dd^{c}|F|^{2})^{p-1}=$$
$$
t^{-2p}\int_{A_{t_{0}}}(dd^{c}|F|^{2})^{p}=
k\,t^{-2p}\int_{|w|<t_{0}}(dd^{c}|w|^{2})^{p}>0
$$
because the function $|w|^{2}$ is strictly plurisubharmonic.
It follows that the map
$\,F:\GA_{j}\to\CC^{p}\sm 0\,$ is not homotopic to a constant.
Then
the same is true for the map $\,F/|F|:\GA_{j}\to\SA^{2p-1}$ to the
unit sphere in $\,\CC^{p}$. If $\,F/|F|$ would be homotopic to a
constant on $\,E\,$ then it would be homotopic to a constant map in
a neighborhood of $\,E\,$ and so on $\,\GA_{j}$ for $\,j\,$ large
enough.
But it is not the case by the proving above, hence
$\,F/|F|:E\to\SA^{2p-1}$ is not homotopic to a constant.
The map
$\,F$ is constant on $\,E'$ and thus it induces the continuous map
$\,h:(E\sm E')\sqcup\circ\to\SA^{2p-1}$, $h(z):=(F/|F|)(z)$ for
$\,z\in E\sm E'$ and $\,h(\circ):=(F/|F|)(E')$ which is also not
homotopic to a constant map (into $\,\SA^{2p-1}$).
And it follows
that $\,H^{2p-1}((E\sm E')\sqcup\circ,\ZZ)\not=0$ (see the same
theorems in \cite{HW}). \vskip1mm

If $\,A\not\sub\widehat{E'}$ we can argue as follows.
By the definition of polynomially convex hull for given
$\,a\in A\sm\widehat{E'}\,$
there is a polynomial $\,f_1\,$ such that $\,f_1(a)=0\,$ and
$\,Re\,f_1>r>0\,$ on $\,E'$.
As $\,\HC^{2p}(E\sm E')=0\,$ there is a polynomial map
$\,f:\CC^n_z\to\CC^p_w\,,\ f=(f_1,...,f_p),$ such that $\,0\in f(A)$
but $\,Re\,f_1>0\,$ on $\,E$.
The rest is the same as in the case $\,\HC^{2p-1}(E')=0\,$ because
the ray $\{Im\,w_1=0,\,Re\,w_1<0,\,w_j=0,\,j=2,...,n\}$ does not
intersect $\,f(E')$.
$\trr$ \vskip3mm

{\bf Lemma 4.} {\it Let $\,E\,$ be a closed subset of a Riemann
surface $\,S\,$ such that $\,H^1(E\sqcup\circ,\ZZ)=0$.
Then the complement $\,S\sm E\,$ is connected.} \vskip2mm

$\trl$ Let $\,a\not= b\in S\sm E$.
Then there is a meromorphic
function $\,f\,$ on $\,S\,$ with only simple zero at $\,a\,$ and
only simple pole at $\,b$.
Let $\,\ga\,$ be a smooth Jordan arc in $\,S\,$
with endpoints $\,a, b$.
Then the multivalued function $\,Log\,f\,$
has a singlevalued holomorphic brunch $\,\log f=\log |f|+i\,arg\,f\,$
in $\,U\sm\ga\,$ for some neighborhood $\,U\supset\ga\,$ homeomorphic
to a disk. \vskip1mm

Let $\,\rho\,$ be a smooth function on $\,S\,$ with zero-set
$\,\ga\,$ such that $\,0\le\rho< 1\,$ and $\,\rho(z_{n})\to 1\,$
for any sequence $\{z_{n}\}\sub S\,$ without cluster points.
Then
there is $\,r>0\,$ such that $\{\rho\le r\}\sub U$.
Let $\,\la(t)$
be a smooth function on $\,\RR_{+}$ such that $\,\la(t)=1\,$ for
$\,0\le t\le r/2\,$ and $\,\la(t)=0\,$ for $\,t\ge r$.
Define
$\,f_{1}:=\exp((\la\circ\rho)\log\,f)$ on $\,S\sm\ga,\,f_{1}=f\,$
on $\,\ga$. Then $\,f_{1}=f\,$ in $\{\rho\le r/2\}$ and
$\,f_{1}=1\,$ in $\{\rho\ge r\}$.
Extend $\,f_{1}$ onto
$\,E\sqcup\circ\,$ setting $\,f_{1}(\circ):=1$.
Then $\,f_{1}$ is
continuous and zero-free on $\,E\sqcup\circ\,$.
By Bruschlinsky
theorem it has continuous logarithm $\,\log f_{1}$ on
$\,E\sqcup\circ\,$ such that $(\log f)(\circ)=0$.
By continuity
$\,\log f=0\,$ in a neighborhood $\,V\ni\circ\,$ in
$\,E\sqcup\circ$. Let $\,V\sub\{\rho>r\}$ be an open subset of
$\,S\,$ such that $\,V\cap E\sub V_{0}$.
Then $\,\log f_{1}$
extended by zero on $\,V\sm E\,$ is continuous logarithm of
$\,f_{1}$ on $\,V\cup E$.

Fix a continuous complete distance {\sl dist} on $\,S\,$ and
denote by $\,\omega\,$ the
modulus of continuity of $\,f_{1}$ on $\,S\sm V$.
If $\,z\in S\,$
is such that $\,\omega(dist(z,E))<(1/4)\min_{E\sm V}|f_{1}|$ then
we define
$$(\log f_{1})(z):=(\log f_{1})(z_{0})+
\log\left(1+\frac{f_{1}(z)-f_1(z_{0})}{f_1(z_{0})}\right)$$
where
$\,z_{0}$ is a nearest point to $\,z\,$ on $\,E\,$ and $\,\log
(1+\eta)$ in $\{\eta\in\CC: Re\,\eta>0\}$ is the continuous brunch
of $\,Log\,\eta\,$ defined by the condition $\,\log 1=0$.
(It
follows from the definition that $\,\log f_{1}(z)$ does not depend on
the choice of nearest point in $\,E$.)
Thus we have defined a continuous logarithm of $\,f_{1}$ in a
neighborhood $\,V_{1}\supset V\cup E$.\vskip1mm

Now assume that $\,E\,$ divides $\,a\,$ and $\,b\,$ and denote by
$\,W\,$ the connected component of $\,S\sm E\,$ containing $\,a$.
Let $\,\rho_{1}$ be a smooth function on $\,W\,$ with zero set
$\{a\}$ such that $\,0\le\rho_{1}<1\,$ and $\,\rho_{1}(z_{n})\to
1\,$ for any sequence $\{z_{n}\}\sub W\,$ without cluster points in
$\,W$.
Choose $\,r_{1}>0\,$ such that $\{\rho_{1}\le
r_{1}\}\sub\{\rho<r/2\}$, then $\,r_{2}\in (r_{1},1)$ such that
$\{\rho_{1}\ge r_{2}\}\sub V_{1}$ and the levels
$\{\ga_{j}:\rho_{1}=r_{j}\}$ are smooth.
The form
$\,f_{1}^{-1}df_{1}$ is closed in $\,W\sm a\,$ hence
$$ 2\pi i=\int_{\ga_{1}} f^{-1} df=\int_{\ga_{1}}f_{1}^{-1} df_{1}=
\int_{\ga_{2}}f_{1}^{-1} df_{1}=\int_{\ga_{2}} d\log f_{1}=0$$ and
this contradiction proves that $\,S\sm E\,$ indeed is connected.
$\trr$\vskip3mm

{\bf Corollary 2.} {\it Let $\,A\,$ be irreducible one-dimensional
analytic set in $\,\CC^{n}$ and $\,\Sigma\,$ is relatively
closed subset of
$\,A\,$ such that $\,\Sigma\sqcup\circ\,$ is simply co-connected.
Then the complement $\,A\sm\Sigma$ is connected.}\vskip2mm

$\trl$ Let $\,\pi:S\to A\,$ be the normalization of $\,A\,$ i.e.
$S\,$ is Riemann surface and $\,\pi\,$ is proper holomorphic map
one-to-one over $\,reg\,A\,$ and such that \#$\,\pi^{-1}(a)$ for
any $\,a\in A\,$ equals to the number of irreducible germs of
$\,A\,$ at the point $\,a$.
Then $\,\pi\,$ extends to
continuous map of compactifications $\,S\sqcup\circ'\to
A\sqcup\circ\supset\Sigma\sqcup\circ$.
Set $\,E:=\pi^{-1}(\Sigma)$
and show that $\,E\sqcup\circ'$ is simply co-connected.
If
$\,\Sigma\cap\,sng\,A=\varnothing\,$ then there is nothing to
prove.
In general $\,\Sigma\cap\,sng\,A=\{a_{1},a_{2},...\}$ is
discrete set and the only possible cluster point of this set in
$\,\Sigma\sqcup\circ\,$ is $\,\circ$.
Similarly,
$E\cap\pi^{-1}(sng\,A)$ is discrete with only possible cluster
point $\,\circ'$ in the compact set $\,E\sqcup\circ'$.

Let $\,f\,$ be a continuous function without zeros on
$\,E\sqcup\circ',\,f(\circ')=1$, extended by continuity onto
$\,S\sqcup\circ'$.
Then $\,f\,$ does not vanish in a neighborhood
$\,V\supset E\sqcup\circ'$ in $\,S\sqcup\circ'$.
Let $\,U_j\ni
a_{j}$ be mutually disjoint neighborhoods of $\,a_{j}$ in
$\,\CC^{n}$ such that $\,\pi^{-1}(U_j\cap A)$ are disjoint unions
of holomorphic disks $V_{jk}\Subset V$ and the continuous variation
of argument of $\,f\,$ on each disk $\,V_{jk}$ is less than $\,\pi$.
Let $\,\la\in C^{\infty}_{0}(\cup U_{j}),\,0\le\la\le
1,\,\la(a_{j})=1\,$ and $(\log f)_{jk}$ are continuous brunches of
logarithm in $V_{jk}$ such that $|Im\,(\log f)_{jk}|<2\pi$.
Then
$\,f_{0}:=f\circ\pi^{-1}$ on $(A\sqcup\circ)\sm(\cup U_{j})$,
$\,f_{0}:=\exp((1-\la))(\log f)_{jk}\circ\pi^{-1}$ in
$\,\pi(V_{jk})$ is continuous function without zeros on
$\,\Sigma\sqcup\circ,\,f_{0}(a_{j})=1$.
As $\,\Sigma\sqcup\circ$ is
simply co-connected there is continuous logarithm $\,\log f_{0}$ in
a neighborhood of $\,\Sigma\sqcup\circ\,$ in $\,A\sqcup\circ$.
Set $\,h:=f/(f_{0}\circ\pi)$.
Then $\,h\,$ is continuous on
$\,S\sqcup\circ'$, equals to 1 on $(S\sqcup\circ')\sm(\cup V_{jk})$
and the variation of argument of $\,h\,$ on each $V_{jk}$ is less
than $\,\pi$.
It follows that $\,h\,$ has continuous logarithm
$\,\log h\,$ on $\,E\sqcup\circ'$ vanishing on
$(E\sqcup\circ')\sm(\cup V_{jk})$.
Thus $\,f\,$ has on $\,E\sqcup\circ'$ continuous logarithm equal
to $(\log f_{0})\circ\pi+\log h$.
By Bruschlinsky theorem $\,E\sqcup\circ'$
is simply co-connected, by Lemma 4  $\,S\sm E\,$ is connected hence
$\,A\sm\Sigma=\pi(S\sm E)$ is connected too.$\trr$\vskip4mm

{\bf 3. The proof of Theorem 1.} We prove more general statement.
\vskip2mm

{\bf Proposition 1.} {\it Let $\,\Sigma\,$ be a bounded relatively
closed subset of a domain $\,D\sub\CC^n$ and
$\,\Sigma_b:=\bar\Sigma\cap\pa D$.
Assume that

$1)\ \HC^2(\Sigma)=0\,,$

$2)\ H^1(\Sigma\sqcup\circ,\ZZ)=0\,$ and

$3)\ \hat\Sigma_b^r\cap\Sigma=\varnothing$.

\np Then $\,\Sigma\,$ is $1$-removable.}
\vskip2mm

Here $\,\hat X^r$ means the rationally convex hull of a compact set
$\,X\sub\CC^n$ that is the set $\{z:P(z)\in P(X)\ \text{for any
polynomial}\ P\,\}$.
Equivalent definition: $\hat X^r$ is the set
of points $\,z\in\CC^{n}$ such that $\,|r(z)|\le max_X\,|r|\,$ for
every rational function $\,r\,$ with poles outside of
$\,X\cup\{z\}$.
Compact sets of zero area ($\,\HC^2$) in $\,\CC^n$ are rationally
convex (coincide with hulls) and thus Theorem 1 follows from
Proposition 1. \vskip2mm

$\trl$
We can assume that $\,D\,$ is bounded.
Let $\,A\sub D\sm\Sigma\,$
be purely 1-dimensional relatively closed analytic subset.
By
Corollary 1 we can assume that $\,\Sigma=(\bar A\sm A)\cap D$.

Let $\,\rho,\ 0\le\rho<1,$ be a smooth function in $\,D$ equal to 0
on $\,\SI\,$ and tending to 1 as $\,z\to\pa D\sm\bar\SI$.
By Sard
theorem almost all levels $\{\rho=t\}$ are smooth hypersurfaces.
The set of singular points of $\,A\,$ is discrete, hence almost all
levels $\,A\cap\{\rho=t\}$ are smooth one-dimensional manifolds.
(They can not be empty because $\,\SI\sub\bar A\,$ and otherwise
the boundary of $\,A\,$ would be
contained in $\,\bar\Sigma\,$ in contradiction with Lemma 3 and the
condition 2.)
Fix such a $\,t\in(0,1)$ with these two properties
and set $\,\Omega:=\{z\in D\sm\hat\Sigma_b^r:\rho<t\},\
\ga:=\pa\Omega\cap A$.\vskip1mm

Now we construct one-dimensional relatively closed analytic subset
$\,A^0\sub\Omega\,$ containing $\,A\cap\Omega$ (the main part of
the proof).\vskip1mm

$\bullet\,$ Let $\,U\,$ be a neighborhood of $\hat\Sigma_b^r$
such that $\,\Sigma\not\sub U$ (see condition 3).
Then there is a compact rational polyhedron
$\,V=\{z\in\CC^n: |p_j(z)|\le 1,\,|q_j(z)|\ge 1,\,j=1,...,N\}$ with
polynomials $\,p_j, q_j$ such that $\hat{\Sigma_b}^r\sub V\sub U$.
Let $\,q(z)\,$ be a polynomial dividing by $\,q_1,...,q_n$ and such
that $\{q=0\}\cap\Sigma=\varnothing$ (it exists due to condition 1).
As $\,q(\Sigma\cup\ga)\sub\CC\,$ has zero area and
$\{q=0\}\cap\Sigma=\varnothing\,$ we can assume (slightly varying
$\,q_j,\,q\,$ and $\,V$) that
$\{q=0\}\cap(\bar\Sigma\cup\ga)=\varnothing$.
Choose $\,\eps>0\,$ so small that
$\{|q|\le\eps\}\cap(\bar\Sigma\cup\ga)=\varnothing$ and set
$\,\Omega':=\Omega\cap\{|q|>\eps\},\,A'=A\cap\Omega'$ and
$\,V'=V\cap\{|q|>\eps\}$.
For $\,\eps\,$ small $\,A'$ is obtained from $\,A\cap\Omega\,$ by
removing of finite number of closed holomorphic disks.

Let $\,M\,$ be the hypersurface $\,w\cdot q(z)=1\,$ in
$\,\CC^{n+1}=\CC^n_z\times\CC_w\,$.
Denote by $\,\pi\,$ the projection $(z,w)\mapsto z\,$ and lift the
picture in $\,\CC^n\sm\{q=0\}$ onto $\,M$ setting $\,\tilde X\,$
be the subset of $\,M\,$ with given projection $\,X$.
Then
$\,\widetilde{V'}=M\cap\{|p_j(z)|<1,\,|w\cdot(q/q_j)|\le 1,\,j=1,...,N\}$
is compact polynomially convex subset of $\,\CC^{n+1}$ and
$\,\tilde\ga\,$ is smooth one-dimensional manifold closed in
$\,\CC^{n+1}\sm\widetilde{\Sigma_b}$.
Show that $\widetilde{A'}$ is contained in polynomially convex hull
of the compact set $\,Y:=\tilde\Sigma_b\cup\widetilde{\ga'}$.

Assume it is not so.
Then there is a point $\,a\in\widetilde{A'}$ and a polynomial $\,P\,$
in $(z,w)$ such that $\,P(a)=1,\,|P|<1\,$ on $\,Y\,$ and the set
$\,P(\widetilde{A'})$ contains a neighborhood of $1$ in $\,\CC$.
As $\,P(\tilde\Sigma)$ has zero area there is $\,a'\in\widetilde{A'}$
such that $|P(a')|>1$, in particular, $a'\not\in\hat Y$.
The boundary of $\,\widetilde{A'}$ is contained in $\,Y\cup\tilde\Sigma$.
As $\,\widetilde{A'}\ni a'\not\sub\hat Y\,$ we obtain by Lemma 3 (with
$\,E'=Y,\,E=Y\cup\tilde\Sigma$) that
$\,H^1(\tilde\Sigma\sqcup\circ,\ZZ)\not=0\,$ in contradiction with
the condition 2 of Proposition 1.
Thus $\,\widetilde{A'}\sub\hat Y$.

The set $\,Y\,$ is contained in $\widetilde{V'}\cup\widetilde{\ga'}=:X$
and $\,\widetilde{V'}$ is polynomially convex.
By Stolzenberg theorem \cite{Sb} $\,\hat X\sm X=:\widetilde{A''}$ being
non-empty (it contains $\widetilde{A'}$) is
bounded purely one-dimensional analytic set with boundary in $\,X$.
As $\,\hat Y\sub\hat X\,$ this analytic set contains
$\,\tilde{A'}\sm X$.
Denote by $\,A''$ the projection of $\widetilde{A''}$ into $\,\CC^n_z$.
Then $\,A''\cup\,(A\cap\Omega\cap\{|q|\le\eps\})$ is relatively
closed analytic subset of $\,\CC^n\sm(V\cup\ga)$ containing
$(A\cap\Omega)\sm V$.
(It is obtained from $\,A''$ by adding of finite number of holomorphic
disks which were removed from $\,A\cap\Omega\,$ before the lifting
onto $\,M$.)
Denote by $\,A_V$ the union of all irreducible components of this set
having relatively open parts of intersections with $\,A\cap\Omega\sm V$.
By the uniqueness theorem for analytic sets (Proposition 5.6.1
in \cite{C}) the set $\,A_V$ does not depend on the choice of $\,q\,$
and $\,\eps\,$ with properties pointed above.

Now we can represent $\,\hat\Sigma^r_b$ as the intersection of decreasing
sequence of rational polyhedrons
$\,V_k=\{z\in\CC^n: |p_{kj}(z)|\le 1,\,|q_{kj}(z)|\ge 1,\,j=1,...,N_k\}$
such that $\,q_k:=\prod_j q_{kj}$ divides $\,q_m$ if $\,m>k$.
By the construction above we have purely one-dimensional analytic sets
in $\,\Omega\sm(V_k\cup\ga)$ containing $(A\cap\Omega)\sm V_k\,.$
Let now $\,\tilde V_k$ be the lifting of $\,V_k$ onto
$\,M_m:w\cdot q_m(z)=1,\,m>k,$ and $\,\eps_m>0\,$ is chosen as above.
As $\,q_k$ divides $\,q_m$ the set $\,\tilde V_k\cap\{|w|\le 1/\eps_m\}$
is polynomially convex and contains $\,\tilde V_m\,$, the lifting of
$\,V_m$.
Thus the polynomially convex hull of
$(\tilde V_k\cap\{|w|\ge 1/\eps_m\})\cup\tilde\ga'_m$ contains the hull
of $\,\tilde V_m\cup\tilde\ga'_m$ and it follows that $\,A_k:=A_{V_k}$
contains $\,A_m\sm V_k$.
Both these sets contain $(A\cap\Omega)\sm V_k$ and thus they coincide
by the uniqueness theorem.
It follows that the union $\,\cup_k A_k=:A^0$ is purely one-dimensional
analytic set relatively closed in $\,\Omega\,$ and
containing $\,A\cap\Omega$. $\bullet$
\vskip1mm

By Lemma 1 we can assume that $\,A\cap\Omega\,$ is irreducible.
Let $\,A^1$ be the irreducible component of $\,A^0\cap\Omega$
containing $\,A\cap\Omega$.
Then $\,\Sigma\sub A^1$ because $\,\Sigma\sub\bar A$.
By Lemma 4 the set $(reg\,A^1)\sm\Sigma\,$ is connected.
But then it coincides with its open subset $\,A\cap\Omega\sm sng\,A^1$
because $\,A\cap\Omega\,$ is closed in $\,\Omega\sm\Sigma$.
It follows that $\,\bar A\cap\Omega=A^1$.
$\trr$
\vskip3mm

Theorem 2 for $\,p=1\,$ follows from Proposition 1 because
$\,\hat\Sigma^r_b\sub\widehat{\Sigma_b},\ D\,$ is
arbitrary and we can substitute $\,D\,$ onto (connected components of)
$\,D\sm\widehat{\Sigma_b}\,$.
\vskip3mm

{\sf Remark.} If one knows that the hulls of $\,V_k\cup\ga\,$
are contained in $\,\Omega\,$ then the proof of the last
part does not need in Corollary 2.
Indeed, if then $\,A'$ is an irreducible component of
$\,A^0\sm\Sigma\,$
then by Lemma 3 it has non-empty open part of its boundary placed
on $\,\ga'$.
By boundary uniqueness theorem as in the proof
above $\,A'$ has non-empty open intersection with $\,A\cap\Omega$.
As $\,A'$ and $\,A\cap\Omega\,$ are closed in $\,\Omega\sm\Sigma\,$
and $\,A'$ is irreducible it
follows that $\,A'\sub A$.

\np The property that the hull is contained in $\,\Omega\,$
(what is used implicitly in the
proof of Theorem 3.8.18 in \cite{S}) is fulfilled, say, for
polynomially
convex (Runge) domains $\,D\,$ but it is not valid in general,
$\,A''\not\sub D\,$ and $\,A''\cap\pa\Omega\not\sub\ga\,$ for
common domains $\,D\,$ because $\,\ga\,$ can be not connected
even if $\,A\cap\Omega\,$ is irreducible (purely one-dimensional
specificity).
In the proof
above this difficulty is overcame due to Lemma 4 and Corollary 2.
\vskip4mm

{\bf 4. Examples.} No one of essential conditions of Theorem 1 is
necessary for the removability of $\,\SI$. \vskip1mm

\np{\sf 1.} The circle $\,\ga=\{z_2=0,\,|z_1|=1,\, y_1\le 0\}\cup
\{z_2=z_1^2-1,\,|z_1|=1,\,y_1\ge 0\}\,$ in $\,\CC^{2}$ is not
simply co-connected but it is removable for purely one-dimensional
analytic sets in $\,D:|z|<3$.
Indeed, let $\,A$ be such a set,
irreducible and closed in $\,D\sm\ga$.
As the singularities of zero
length are removable then  $\,\pa A\,$ contains a part of
$\,\ga\,$ of positive length and by the boundary uniqueness theorem
$\,A\,$ coincides either with $(D\sm\ga)\cap\{z_2=0\}$ and then
$\,\bar A\cap D=\{z_2=0\}\cap D\,$
or with $(D\sm\ga)\cap\{z_2=z_1^2-1\}$ and then
$\,\bar A\cap D=\{z_2=z_1^2-1\}\cap D$. \vskip1mm

\np{\sf 2.} Let $\,E\,$ be a closed totally disconnected set of
finite length in the unit disk $\,\DD$ and
$\,\SI:=E^n\sub\DD^n$.
Let $\,A\,$ be an irreducible
one-dimensional relatively closed analytic subset of
$\,\DD^n\sm\SI\,$ and $\,a\in\SI\,$ is its cluster point.
One can assume that $\,A\,$ is
not contained in a plane $\,z_1=c_1$ for $\,c_1\in\DD$.
Then
$\,A\cap\{z_1=c_1\}$ is a discrete set and its union with
$\,\SI\cap\{z_1=c_1\}$ is closed and totally disconnected.
Thus
there exists a neighborhood $\,V\,$ of the point $(a_2,...,a_n)$ in
$\,\DD^{n-1}$ such that $\,a_1\times\pa V\,$ does not intersect
$\,(A\cup\SI)\cap\{z_1=c_1\}$. Let $\,r>0\,$ be so small that
$\,U:=\{|z_1-a_1|<r\}\Subset\DD\,$ and $\,U\times\pa V\,$ does not
intersect $\,A\cup\SI\,$ also.
Then the projection $\,z\mapsto z_1$
of the set $(A\cup\SI)\cap(U\times V)$ is proper, hence there
exists $\,k\in\NN\,$ such that the number of points in
$\,A\cap(c_1\times V)$ counting with multiplicities is equal to
$\,k\,$ for all $\,c_1\in U\sm E$.
As $(A\cup\SI)\cap(c_1\times V)$ is totally disconnected for each
$\,c_1\in U$ then \#$\,\bar A\cap(c_1\times V)\le k,\,c_1\in U$,
and at each point $\,c\in\bar A\cap(U\times V)$ the multiplicity of
$\,z_1$ is well defined so that the number of points in the
intersection of $\,\bar A\cap(U\times V)$ and $\{z_1=c_1\}$
counting the multiplicity is equal to $\,k$.
Thus the projections
of the set $\,\bar A\cap(U\times V)$ into $\,\CC^{\,2}_{z_1
z_j}\,,\ j=2,...,n,$ are given by corresponding equations
$\,z_j^k+s_{j1}(z_1)z_j^{k-1}+\cdots+s_{jk}(z_1)=0\,,$ where the
functions $\,s_{ji}$ are continuous in $\,U\,$ and holomorphic in
$\,U\sm E$.
As the length of $\,E\cap U\,$ is finite the functions
$\,s_{ji}$ are holomorphic in $\,U\,$ (see e.g. Theorem A$\,$1.5 in
\cite{C}) hence, the projections $\,\bar A\cap(U\times V)$ to
$\,\CC^{\,2}_{z_1 z_j}$ are analytic.
It follows evidently the
analyticity of $\,\bar A\cap(U\times V)$ and $\,\bar A\cap\DD^n$.
In this example $\,n>1\,$ is arbitrary and $\,\SI\,$ is removable
in spite of $\,\HC^n(\SI)>0$. \vskip3mm

{\bf 5. The proof of Theorem 2.} We need in some properties of
solutions of the Plateau problem for analytic sets (see \S 19.3 in
\cite{C}).

A real $C^{1}$-manifold $\,\La\,$ of dimension $2p-1,\,p>1,$ in a
complex manifold $\,M\,$ is called maximally complex if the
dimension of its complex tangent space $T_{a}\La\cap i\,T_a\La\sub
T_{a}M\,$ has complex dimension $p-1$ at every point $\,a\in\La$.

A closed subset $\,\GA\,$ of a complex manifold $\,M\,$ is called
{\sf maximally complex cycle} (of dimension $2p-1,\,p>1$) if the
measure
$\,\HC^{2p-1}|_{\GA}$ is locally finite and there is a closed
(maybe empty) subset $\,\sigma\sub\GA\,$ of zero
$\,\HC^{2p-1}$-measure such that $\,\GA\sm\sigma\,$ is a smooth
($C^1$) oriented maximally complex manifold of dimension $2p-1$
and the current of integration on $\,\GA$ (of smooth differential
forms of degree $2p-1$ with compact supports in $\,M$) is closed.
Such a
cycle is called {\sf irreducible} if it contains no proper maximally
complex cycle of the same dimension.

If $\,A\,$ is a closed purely $\,p$-dimensional analytic subset in
$\,M\,$ and $\,\rho\,$ is a real smooth function on $\,M\,$ then
for almost every $\,t\in\rho(A)$ the set
$\,\GA_{t}:=A\cap\{\rho=t\}$ with smooth part oriented as the
boundary of $\,A_{t}:=A\cap\{\rho<t\}$ (with canonical orientation
corresponding to the complex structure on $\,M$) is a maximally
complex cycle of dimension $2p-1$ (Proposition 14.3.1 in \cite{C}).
In this case Stokes formula is valid:
$\,\int_{A_{t}}d\phi=\int_{\GA_{t}}\phi\,$ for every smooth form of
degree $2p-1$ on $\,M$ (Theorem 14.3 in \cite{C}).
As the complex
dimension of the set of singular points in $\,A\,$ is not more than
$p-1$ then for almost every $\,t\in\rho(A)$ the singular set
$\,\sigma_{t}$ from the definition has locally finite
$\,\HC^{2p-2}$-measure and the irreducible components of $\,\GA_{t}$
are precisely the closures of connected components of
$\,\GA_{t}\sm\sigma_{t}$.

If $\,Y\,$ is polynomially convex compact subset in $\,\CC^{n}$ and
$\,\GA\,$ is a bounded maximally complex cycle in $\,\CC^{n}\sm Y$
of dimension $2p-1,\,p>1,$ then by generalized Harvey -- Lawson
theorem (Theorem 19.6.2 in \cite{C}) there exists a bounded closed
in $\,\CC^{n}\sm(Y\cup\GA)$ purely $\,p$-dimensional analytic set
$\,A'$ such that $\,\GA\sub\overline{A'}$.
By the theorem on
boundary regularity (Theorem 19.1 in \cite{C}) and boundary
uniqueness theorem (Proposition 19.2.1 in \cite{C}) the set $\,A'$
is irreducible if such is the cycle $\,\GA$.

Reformulate Theorem 2 for $\,p>1\,$ taking in mind Lemma 1 (the case
$\,p=1\,$ is already considered in Proposition 1).
\vskip2mm

{\bf Proposition 2.} {\it Let $\,D\,$ be a domain in $\,\CC^{n},\
\SI\sub D\,$ is bounded relatively closed subset,
$\SI_{b}:=\bar\SI\cap\pa D\,$ and
$\,A\,$ is relatively closed irreducible analytic set of dimension
$\,p>1\,$ in $\,D\sm\SI$.
Assume that
\vskip1mm

$1)\ \,\HC^{2p}(\SI)=0\,,$

$2)\ \,H^{2p-1}(\SI\sqcup\circ,\ZZ)=0\,$ and

$3)\ \,\widehat{\SI_b}\cap\SI=\varnothing\,$.
\vskip1mm

\np Then $\,\bar A\cap D\,$ is analytic and $\,\SI\,$ is $\,p$-removable.}

\vskip2mm

$\trl$ Substituting $\,\SI\,$ onto $\,\SI\cap\bar A\,$ one can
assume (due to Corollary 1) that $\,\SI\sub\bar A\,$ and we
suppose this in what follows.

Let $\,\rho,\ 0\le\rho<1,$ be a smooth function in $\,D$ equal to 0
on $\,\SI\cup(\widehat{\SI_b}\cap D)\,$ and tending to 1 as
$\,z\to\pa D\sm\widehat{\SI_b}$.
By Sard
theorem for analytic sets (Proposition 14.3.1 in \cite{C}) almost
all levels $\,A\cap\{\rho=t\}$ are either empty or maximally
complex cycles in $\,D\sm\widehat{\SI_b}\,$ of dimension $2p-1$.
If $\,A\cap\{\rho=t\}$ is empty then either
$\,A\cap\{\rho<t\}$ is empty or it is nonempty closed in $\,A$.
In the
first case there is nothing to prove ($\,\SI\cap\bar A=\varnothing$)
and in the
second case $\,A\,$ is contained in $\{\rho<t\}$ because $\,A\,$ is
irreducible.
But then $\,\pa A\sub\bar\SI\,$ in contradiction with
the condition 2 and Lemma 3.
Thus one can assume that almost all levels
$\,A\cap\{\rho=t\}$ are maximally complex cycles in $\,D\,$ of
dimension $\,2p-1$.

From Lemma 3 with $\,E'=\widehat{\SI_b}$ and $\,E=E'\sqcup\SI\,$
we obtain that $\,\pa A\not\sub \widehat{\SI_b}\cup\SI,$ in
particular, $A\not\sub\widehat{\SI_b}$.
Hence, there is $\,t_1>0\,$ such that
$\,A^t:=A\cap\{\rho<t\}\not\sub\widehat{\SI_b}$ and
$\,\pa A\not\sub\widehat{\SI_b}\cup\SI$ for $\,t\ge t_1\,$.
Fix $\,t_0>t_1$ such that $\,A^{t_0}$ is analytic set with
maximally complex cycle-border
$(\pa A^{t_0})\sm(\widehat{\SI_b}\cup\SI)$.
Denote by $\,\GA_0$ an irreducible component of this cycle and
by $\,A_0$ the irreducible component of $\,A^{t_0}\sm\widehat{\SI_b}$
whose boundary contains $\,\GA_0\,$.
Set $\,X:=\widehat{\SI_b}\cup\GA_0\,$.

By the generalized Harvey -- Lawson theorem (Theorem 19.6.2 in
\cite{C}) there exists bounded irreducible $p$-dimensional
analytic set
$\,A'$ in $\,\hat X\sm X\,$ with boundary in $\,X\,$ such that
$\,\GA_0\sub\pa A'$.
By the boundary uniqueness theorem (Proposition
19.2.1 in \cite{C}) and Shiffman theorem
there is a neighborhood
$\,U\supset\GA_{0}\,$ in $\,D\sm\SI\,$ and a relatively closed set
$\,\sigma\sub\GA_{0}\,$ (maybe empty) of zero
$\,\HC^{2p-1}$-measure such that either
$(A'\cup\GA_{0}\cup A_0)\cap(U\sm\sigma)$ is analytic or
$\,A'\cup\GA_{0}\sm\sigma\,$ is a smooth manifold with boundary
$\,\GA_0\sm\sigma\,$ in $\,U\sm\sigma$.

In the first case $\,A'\cup\GA_0\cup A_0$ is $\,p$-dimensional
analytic set with boundary in $\,\widehat{\SI_b}\cup\SI\,$ what
is impossible by Lemma 3 and the condition 2.
In the second case $\,A_0\cap A'$ have $\,p$-dimensional
intersection what follows that $\,A_0\sub A'$ because $\,A_0$ is
irreducible and $\,A'$ is closed in $\,\CC^n\sm X$.
By the boundary uniqueness theorem $\,A'$ and $\,A_0$ coincide in
a neighborhood of $\,\GA_0\sm\sigma$.
It follows that $\,A'\sm(A_0\cup(\SI\cap\overline{A_0}))$
if it is non-empty is $\,p$-dimensional analytic set with boundary
in $\,\widehat{\SI_b}\cup\SI\,$ what again is impossible.
Hence, $A'=A_0\cup\,(\SI\cap\overline{A_0})$ and thus
$\,\overline{A_0}$ is analytic in a neighborhood of
$\,\SI\,$ in $\,D$.
$\trr$
\vskip3mm

{\sf Remark.} Proposition 2 can be generalized by weakening the
last condition to $\,\hat\SI^r_b\cap\SI=\varnothing\,$ (lifting
the picture into $\,\CC^{n+1}$ as in the proof of Proposition 1).
But in case $\,p>1\,$ more natural would be not the rational
convexity but some convexity with respect to polynomial mappings
to $\,\CC^p$.
But the method used above does not work for such
more weak conditions.

\np The proof of Proposition 2 is essentially simpler than that one
for Proposition 1 because of two crucial advantages of the case
$\,p>1$.
The first one is in that then one does not need in the proof of the
{\sf existence} of analytic set $\,A'$ with boundary in
$\,\widehat{\SI_b}\cup\GA_0$ because it follows from generalized
Harvey -- Lawson theorem.
The second one is in that the cicle $\,\GA_0$ is irreducible if
(and only if) $A'$ is irreducible.
\vskip4mm

{\bf 6. Corollaries.} Several consequences in spirit of \S$\,$3.8
in \cite{S}. First one is about the "infection"$\,$ property of
removability (see Theorem 18.2.1 in \cite{C}). \vskip2mm

{\bf Corollary 3.} {\it Let $\,\ga\,$ be an open Jordan arc
relatively closed in a domain
$\,D\sub\CC^n$ such that $\,\HC^2(\ga)=0$.
Let $\,A\,$ be purely one-dimensional analytic set in
$\,D$ such that $\,(\bar A\sm A)\cap D\sub\ga$.
If $\,\bar A\,$ is
analytic $($maybe empty $)$ in a neighborhood of a point
$\,a\in\ga\,$ then
$\,\bar A\cap D$ is analytic set.} \vskip2mm

$\trl$ Let $\,\ga:(0,1)\to D\,$ be a parametrization, $\ga(0)=a$.
For every $\,\eps\in(0,1/2)$ there is a domain $\,D_\eps\sub D\,$
such that $\,\ga\cap D_\eps=\{\ga(t):\eps<t<1-\eps\}$ and
\#$\,\ga\cap D_\eps=2$.
By the condition there
exists an arc $\,\ga'\sub\ga\cap D_\eps$ such
that $\,\bar A\,$ is analytic in its neighborhood.
The set $(\ga\sm\ga')\cap D_\eps$ if it is non-empty has not more
than two connected components and each of them satisfies the
conditions of Stout theorem (in domain $\,D_\eps$).
It follows that $\,\bar A\cap D_\eps$ is analytic for all $\,\eps$.
$\trr$ \vskip2mm

The corresponding statement for $\,p>1\,$ is true for connected
$\,C^1$-manifolds of dimension $\,2p-1$ (instead of $\,\ga\,$;
see Theorem 18.2.1 in \cite{C}).
For general topological manifolds of zero $\,\HC^{2p}$-measure
it is verisimilar that the statement is valued also but the
proof like given above does not work for $\,p>1\,$ because of
non-checked in this case condition 3 in Proposition 2.
\vskip3mm

{\bf Corollary 4.} {\it Let $\,\SI\,$ be a closed subset of zero
$\,\HC^{2p}$-measure in a domain $\,D\sub\CC^n$ such that
$\,H^{2p-1}(\SI,\ZZ)=0\,$ and each
connected component of $\,\SI\,$ is compact $($say,
$\SI\,$ is totally disconnected$)$.
Then $\,\SI\,$ is $\,p$-removable.}
\vskip2mm

$\trl$ Let $\,A\sub D\,$ be purely $\,p$-dimensional analytic set
such that $\,(\bar A\sm A)\cap D\sub\SI\,$ and $\,\SI'\sub\SI\,$
be a connected component.
Then there exists a neighborhood
$\,U\supset\SI'\,$ with compact closure in $\,D\,$ such that
$\,\pa U\cap\SI=\varnothing$.
Then $\,\SI\cap U\,$ is a compact and $\,H^{2p-1}(\SI\cap U,\ZZ)=0$.
simply co-connected set by Lemma 1.
By Theorem 2 the set $\,\bar A\cap U\,$ is analytic. Thus $\bar
A\cap D\,$ is analytic in a neighborhood of each point of $\,D$.
$\trr$ \vskip2mm

And finally a consequence of pseudoconvexity condition on $\,D$.
\vskip3mm

{\bf Corollary 5.} {\it Let $\,D$ be a bounded domain in $\,\CC^n$
such that $\,\bar D$ is the intersection of a decreasing sequence
of pseudoconvex domains and $\,\pa D=\pa(int\,\bar D)$.
Let $\,A\sub D$ be purely $\,p$-dimensional
analytic set closed in $\,D\,$ and such that $\,\HC^{2p}(\pa A)=0$.
Then $\,H^{2p-1}(E,\ZZ)\not=0\,$ for every simultaneously open and
closed subset $\,E\sub\pa A$.}
\vskip2mm

$\trl$ Assume there is be an open-closed subset $\,E\sub\pa A\,$
with $\,H^{2p-1}(E,\ZZ)=0$.
Then there exists a neighborhood $\,U\supset
E\,$ in $\,\CC^n$ such that $\,\pa U\cap\pa A=\varnothing,\
U\cap\pa A=E\,$ and thus $\,A\cap U\,$ is an analytic set closed in
$\,U\sm E$.
By Theorem 2 $\,\tilde A:=(A\cap U)\cup E\,$ is
analytic set in $\,U$.
By the construction, its boundary $\,\GA\,$
is placed in $\,D\cap\pa U\,$ on the positive distance
$\,\delta>0\,$ from $\,\pa D$.
By the condition there exists a
pseudoconvex domain $\,D'\supset \bar D\,$ such that the distance
of $\,E\,$ to $\,\pa D'$ is less then $\,\delta$.
By the maximum
principle for $\,\tilde A\,$ the set $\,E\,$ is contained in the
hull of $\,\GA\,$ convex with respect to the algebra of functions
holomorphic in $\,D'$.
But the distance of this hull to $\,\pa D'$
can not be less than $\,\delta\,$ due to the pseudoconvexity of
$\,D'$ (see \cite{Sh} Theorem ?).
The contradiction shows that $\,E=\varnothing$. $\trr$
\vskip1mm

In particular, no open subset in $\,\pa A\,$ can be totally
disconnected (see \cite{S} Corollary 3.8.20).
However, there
remains open question (even if $\,D\,$ is a ball and $\,p=1$) if
$\,\pa A\,$
can have a connected component consisting of one point.
\vskip4mm

{\bf 7. Comments and questions.} The graphs of holomorphic functions
are special analytic sets and thus Theorem 1,2 can be applied to
removability of singularities of holomorphic functions
(see \cite{C1,C2}).
But in the case of graphs the conditions on singular sets can be
essentially weakened.
Quote for comparing the main result of \cite{C}:
\vskip1mm

\np {\it Let $\,E\,$ be closed subset of a Riemann surface $\,S\,$
and $\,f\,$ is a meromorphic function on $\,S\sm E\,$ such that
$\,\HC^2(C_f(E))=0\,$ and the cluster set $\,C(f,z)\sub\hat\CC\,$
at any point of $\,E\,$ is connected and has connected complement.
Then $\,f\,$ extends to meromorphic function on $\,S$.}
\vskip1mm

\np(Here $\,C(f,z)$ is the set of cluster values of $\,f(\zeta)$ as
$\,\zeta\to z,\,\zeta\in S\sm E\,$ and
$\,C_f(E):=\cup_{\,z\in E}\, \{z\}\times C(f,z)$, the "graph"$\,$ of
cluster values of $\,f\,$ at the points of $\,E$.)

Note that the statement does not follow from Theorem 1 because there
is no condition on boundary behaviour and global topology of
$\,\SI:=C_f(E)$.
Nevertheless, the proof (based on argument principle too) is simpler
due to simplicity of graphs with respect to general analytic sets.
\vskip2mm

An analogy of theorems on removable singularities for functions
continuous on $\,S\,$ and holomorphic on $\,S\sm E\,$ is the following
(Proposition 19.2.1 in \cite{C}).
\vskip1mm

\np {\it Let $\,M\,$ be a connected $(2p-1)$-dimensional
$\,C^1$-submanifold of a complex manifold $\,\Omega,\,p\ge 1,$ and
$\,A_1\,,\, A_2$ are different irreducible $\,p$-dimensional analytic
sets relatively closed in $\,\Omega\sm M\,$ and such that
$\,M=\bar A_1\cap\bar A_2\,$.
Then $\,A_1\cup M\cup A_2$ is analytic subset of $\,\Omega$.}
\vskip1mm

\np In view of Theorems 1,2 the natural question here is if the
condition $\,M\in C^1$ can be weakened to, say, $M\,$ is connected
and $\,\HC^{2p-1}(M)<\infty$ (or even to $\,\HC^{2p}(M)=0?$).
Maybe some additional topological conditions?
The question is open even if $\,M\,$ is topological manifold of
finite $\,\HC^{2p-1}$-measure.
\vskip2mm

As noted after Corollary 3 the third condition in Proposition 2
(which is not necessary at all) is rather restrictive for applications.
The natural desire arises to substitute it by something like
that in Theorem 1 (say, $\HC^{2p}(\bar\SI)=0$) or similar one.
In any case it would be useful to substitute the condition
$\,\widehat{\SI_b}\cap\SI=\varnothing\,$ by one simpler for
checking.


\vskip4mm

\begin{thebibliography}{99}

\bibitem{B}
N.Bruschlinsky,
\emph{Stetige Abbildungen und Bettische Gruppen der Dimensionszahlen 1 und 3,}
Math. Ann., \textbf{109} (1934), 525--537.

\bibitem{C}
E.M.Chirka,
\emph{Complex analytic sets,}
Kluwer Acad.Publ., Dordrecht, 1989.

\bibitem{C1}
-----,
\emph{On the removable singularities of holomorphic
functions,}
Mat.Sbornik, 2016.

\bibitem{C2}
-----,
\emph{On the $\bar\pa$-problem with $L^2$-estimates
on a Riemann surface,}
Proc. Steklov Inst. Math., \textbf{290} (2015), 264--276.

\bibitem{G}
J.Garnett,
\emph{Analytic capacity and Measure,}
Lecture Notes in Math., \textbf{297}
Springer, Berlin, 1972.

\bibitem{HW}
W.Hurewicz and H.Wallmen,
\emph{Dimension theory,}
Princeton Univ. Press, 1941.

\bibitem{Sh}
B.V.Shabat,
\emph{Introduction to Complex Analysis, Part 2В}
Moscow, Nauka, 1985.

\bibitem{Sf}
B.Shiffman,
\emph{On the removable singularities of analytic sets,}
Michigan Math. J., \textbf{15} (1968), 111--120.

\bibitem{Sb}
G.Stolzenberg,
\emph{Uniform approximation on smooth curves,}
Acta Math., \textbf{115} (1966), 185--198.

\bibitem{S}
E.L.Stout,
\emph{Polynomial convexity,}
Progress in Math., \textbf{261},
Birkh\"auser, Boston, 2007.

\end{thebibliography}
\end{document}